\documentclass[11pt]{amsart}
\usepackage{amssymb}
\usepackage{latexsym}
\usepackage{amsmath}
\usepackage{amsthm}
\theoremstyle{plain}
\newtheorem{lemma}{Lemma}
\newtheorem{theorem}{Theorem}
\theoremstyle{remark}
\newtheorem*{remark}{Remark}
\DeclareMathOperator{\Lie}{Lie}
\DeclareMathOperator{\Span}{Span}
\DeclareMathOperator{\Cl}{Cl}
\DeclareMathOperator*{\limind}{lim\,ind}
\renewcommand{\phi}{\varphi}
\newcommand{\ot}{\otimes}
\newcommand{\BC}{{\mathbb{C}}}
\newcommand{\BN}{{\mathbb{N}}}
\newcommand{\BR}{{\mathbb{R}}}
\newcommand{\BT}{{\mathbb{T}}}
\newcommand{\BZ}{{\mathbb{Z}}}
\newcommand{\CA}{{\mathcal{A}}}

\newcommand{\CH}{{\mathcal{H}}}
\newcommand{\CL}{{\mathcal{L}}}
\newcommand{\CX}{{\mathcal{X}}}

\newcommand{\GA}{{\mathfrak{A}}}
\newcommand{\TGA}{{\widetilde{\mathfrak{A}}}}
\newcommand{\GZ}{{\mathfrak{Z}}}

\begin{document}
%

\title{GRADED LIE ALGEBRAS AND DYNAMICAL SYSTEMS}

\author{A. M. Vershik}
\thanks{\parindent=0pt
St. Petersburg Division,
Steklov Mathematical Institute, RAS,\\
27, Fontanka, St. Petersburg, 191011, Russia.\\
Partially supported by grant RFBR 99-01-00098.\\
This paper contains  the subject of my talk on the
Conference in Twente in December 2000.}%

\maketitle

\section*{Introduction}

We consider a class of infinite-dimensional Lie algebras
which is associated to dynamical systems with invariant
measures. There are two constructions of the algebras -- one
based on the associative cross product algebra which considered
as Lie algebra and then extended with nontrivial scalar two-cocycle;
the second description is the specification of the construction
of the graded Lie algebras with continuum root system in spirit of
the papers of Saveliev-Vershik \cite{VS1,VS2,V} which is a generalization
of the definition of classical Cartan finite-dimensional algebras
as well as Kac--Moody algebras.
In the last paragraph we present the third construction for the
special case of dynamical systems with discrete spectrum.
The first example of such algebras was so called sine-algebras
which was discovered independently in \cite{VS1} and \cite{FFZ} and had been
studied later in \cite{GKL} from point of view Kac--Moody Lie algebras.
In the last paragraph of this paper we also suggest a new examples
of such type algebras appeared from arithmetics: adding of $1$ in
the additive  group $Z_p$ as a transformation of the group  of
$p$-adic integers. The set of positive simple roots in this case
is $Z_p$; Cartan subalgebra is the algebra of continuous functions
on the group $Z_p$ and Weyl group of this Lie algebra contains the
infinite symmetric group. Remarkably this algebra is the inductive
limit of Kac--Moody affine algebras of type $A^1_{p^n}$.

\section{Lie algebra generated by automorphism}

\subsection{Associative algebra $\CA(\CX,T)$}

Let $(\CX,\mu)$ be a separable compactum with Borel probability
measure $\mu$ which is positive for any open set $A\subset\CX$
and $T$ is a measure preserving homeomorphism of $\CX$.

It is old and well-known construction of $W^*$-algebra (von Neumann)
and $C^*$-algebra (Gel'fand) generated by $(\CX,\mu,T)$, see for example
\cite{D,ZM}. Algebraically this is an associative algebra $\CA(\CX,T)$
which is semidirect product of $C(\CX)$ and $C(\BZ)$ with the action
of $\BZ$ on $\CX$ and consequently on $C(\CX)$. As a linear space this
is direct sum
\[
\CA(\CX,T)=\bigoplus_{n\in\BZ}C(\CX)\otimes U^n
\]
where $C(\CX)$ is Banach space of all continuous functions on $\CX$,
and $U=U_T$ is linear operator $(U_Tf)(x)=f(Tx)$, $f\in C(\CX)$.  The
multiplication of the monomials is defined by formula
\[
(\phi\otimes U^n)\cdot(\psi\otimes U^m)
=(\phi\cdot U^n\psi)\otimes U^{n+m},
\qquad n,m\in\BZ,\quad\phi,\psi\in C(\CX).
\]
Involution on $\CA(\CX,T)$ is the following:
\[
(\phi\ot U^n)^*=(U^{-n}\bar\phi)\ot U^{-n}.
\]
Completion of $\CA(\CX,T)$ with respect to the appropriate
$C^*$-norm gives a corresponding $C^*$-algebra.

It is possible to include to this construction a 2-cocycle of the action of
$\BZ$ with values in $C(\CX)$ to obtain another $C^*$-algebra which are
unsplittable extensions (see \cite{ZM, VSh}), but we restrict ourself
to the case of the trivial cocycle.

If we use a measure $\mu$ as a state on $\CA(\CX,T)$ and construct
$*$-representation corresponding to this state then $W^*$-closure of
image of $\CA(\CX,T)$ gives us $W^*$-algebra generated by triple
$(\CX,\mu,T)$.

There are two classical representations of the algebra
$\CA(\CX,T)$ --- Koopmans representation (in $L^2_\mu(\CX)$) and
von~Neumann one in $L^2_{\mu\times m}(\CX\times\BZ)$ ($m$ is Haar
measure on $\BZ$). This area called ``algebraic theory of dynamical
systems'' and there are many papers on this. Extremally popular is so called
rotation algebras (also called as "quantum torus") which is associative
$C^*$-algebra generated by irrational rotation of the unit circle.

{\it We want to point out that there is another remarkable algebraic
object which is associated with dynamical systems---some Lie algebras
which are similar to the classical Cartan Lie algebras and to affine
Kac-Moody algebras.} Some nontrivial central extension included in the
definition plays very imporatant role in the whole theory. Upto now
only  the Lie algebras corresponing to rotation algebras were
considered; it was  discovered independently in \cite{VS1} and \cite{FFZ}
(see also \cite{GKL} were  the shift on d-dimensional torus was considered)
and called by physisits "sine algebras" - all those dynamical systems have a
{\it discrete spectrum},

In whole generality  the Lie algebras generated by an arbitrary
dynamical system with invariant measure was briefly defined in \cite{VS2}
and more systematically in \cite{V}. It is interesting that we started
not form the generaltheory of dynamical systems as in the first
definition below
but from the notions presented in our series of papers with
M.~Saveliev \cite{VS1,VS2} were we had defined so called
``$\BZ$-graded  Lie algebras  with continuous root systems''.
Those algebras which we will discuss here were one of
the type of the examples and a special case of $\BZ$-graded
Lie algebras with general root systems.  Below we will describe
explicitly the modern and detailed version of the construction
of Lie algerbas generated by an arbitrary discrete
dynamical systems with invariant measure and then will give
the link between various definitions. One can hope that
this type of algebras can give a new type of invariants
of the dynamical systems as well as new examples of calssical and
quantum integrable systems.

\subsection{Lie algebras $\widetilde\GA(\CX,T)$}
Most interesting case is the case when $T$ is minimal (= each orbit of
$T$ is dense in $\CX$) and ergodic with respect to measure $\mu$ (=
there are no nonconstant $T$-invariant measurable functions).  We
assume this in further considerations.
It is known that if $T$ is minimal (=each $T$-orbit is dense
in $\CX$) then $C^*$-algebra is simple (= has
no proper two-sided ideals), see \cite{ZM}. Algebra $\CA(\CX,T)$
with brackets
\[
[a,b]=ab-ba
\]
will be denoted as $\Lie\CA(\CX,T)$; it is still $\BZ$-graded Lie
algebra and the brackets of monomials are
\[
[\phi\ot U^n,\psi\ot U^m]=(\phi\cdot U^n\psi-\psi\cdot U^m\phi)\ot U^{n+m}.
\]
This algebra has a center.

\begin{lemma}
The center of $\Lie\CA(\CX,T)$ is the set of constants functions
in zero component subspace:
$\GZ=c\ot U^0$, $c\in\BC$.  The complement linear subspace
\[
\CA_0(\CX,T)=\bigoplus_{n<0}C(\CX)\ot U^n\oplus
C_0(\CX)\ot U^0 \oplus
\bigoplus_{n>0}C(\CX)\ot U^n,
\]
where $C_0(\CX)=\bigl\{\,\phi\in C_(\CX):\int_\CX\phi(x)d\mu=0\,\bigr\}$,
is Lie subalgebra which is isomorphic to quotient $\CA(\CX,T)\big/\GZ$
over center $\GZ$.
\end{lemma}

\begin{remark}
The center is not ideal of associative algebra consequently there
is no ``associative'' analogue of this lemma and $\CA_0(\CX,T)$
is not a subalgebra of $\CA(\CX,T)$.
\end{remark}

Now we define a 2-cocycle on $\CA_0(\CX,T)$ with the scalar values
and one-dimensional central expansion of it.

\begin{lemma}
The following formula defines 2-cocycle on $\CA_0(\CX,T)$:
\[
\alpha(\phi\ot U^n,\psi\ot U^m)
=n\int_{\CX}\phi\cdot U^n\psi\,d\mu\cdot\delta_{n+m},
\]
so
\[
\alpha(\phi\ot U^n,\psi\ot U^m)
=\begin{cases}n\int_{\CX}\phi\cdot U^n\psi\,d\mu & \text{if }m=-n,\\
0 & \text{if }m\ne-n. \end{cases}
\]
\end{lemma}

\begin{proof}
We need to check that $\alpha([x,y],z)+\alpha([y,z],x)+\alpha([z,x],y)=0$.
Let $k+l+n=0$. Then
\begin{multline*}
\alpha\bigl([\phi\ot U^k,\psi\ot U^l],\gamma\ot U^n\bigr)+\dots
=\alpha\bigl((\phi\cdot U^k\psi-\psi\cdot U^l\phi)\ot U^{k+l},
   \gamma\ot U^{-k-l}\bigr)+\dots \\
=(k+l)\int_\CX\bigl(\phi\cdot U^k\psi U^{-n}\cdot U^{-n}\gamma
    -\psi\cdot U^l\phi\cdot U^{-n}\gamma\bigr)d\mu+\dots
=0
\end{multline*}
(Dots mean cyclic permutation of indices, we used here the invariance
of measure
$\mu$ under $T$.)
\end{proof}

\begin{remark}
Cocycle $\alpha$ is not cohomologous to zero because it is easy
to check that $\alpha(x,y)$ can not be represented as $f([x,y])$,
with any linear functional $f$.
of $C(\CX)$.
\end{remark}

Let us identify scalars $c$  which are extensions of $\CA_0(\CX,T)$
with scalars $c\in C(\CX)\ot U^0\subset\CA(\CX,T)$. So we can consider
again linear space $\CA(\CX,T)$ as one dimensional nontrivial extension
of Lie algebra $\CA_0(\CX,T)$.  Denote a new Lie algebra by
$\TGA(\CX,T)$.
So, Lie algebra $\TGA(\CX,T)$ as linear space is the same as $\CA(\CX,T)$
but the brackets in $\TGA(\CX,T)$ differ from the brackets in
$\CA(\CX,T)$:
\begin{equation}\label{eq1}
[\phi\ot U^n,\psi\ot U^m]=(\phi\cdot U^n\psi-\psi\cdot U^m\phi)\ot U^{n+m}
  +\int_{\CX}\phi\cdot U^n \psi\,d\mu \cdot \delta_{n+m}
\end{equation}
It means that the center of $\TGA(\CX,T)$ is again scalars
$\BC\cdot1\subset C(\CX)\ot U^0\subset \Lie \CA(\CX,T)$, but now
subspace $\CA_0(\CX,T)$ is not Lie subalgebra and the central extension is
not trivial.

Lie algebra $\TGA(\CX,T)$ is $\BZ$-graded Lie algebra.  We will
give a new definition of it in a framework of Lie algebras with
continuous root systems.  We will call the subspace of $\TGA(\CX,T)$
which consists of $\Span\{\gamma\ot U^{-1}\}\oplus\Span\{\phi\ot U^0\}
\oplus\Span\{\psi\ot U^1\}$, $\phi,\psi,\gamma\in C(\CX)$,
a ``local subalgebra''.  Here are the brackets for local part of
$\TGA(\CX,T)$ are
\begin{equation}\label{eq2}
\begin{split}
[\phi_1\ot U^0,\phi_2\ot U^0] &{}=0,\\
[\phi\ot U^0,\psi\ot U^{\pm1}]&{}
=\pm\bigl((\phi-U\phi)\cdot\psi\bigr)\ot U^{\pm1}
=\bigl((I-U)\phi\cdot\psi\bigr)\ot U^{\pm1}, \\
[\phi\ot U^{+1},\psi\ot U^{-1}]&{}
=(\phi\cdot U\psi-\psi\cdot U\phi)\ot U^0
  +\int_\CX(\phi\cdot U\psi)d\mu\cdot c.
\end{split}
\end{equation}
The middle term of local algebra ($\{\phi\ot U^0; \phi \in C(\CX) \}$)
is by definition Cartan subalgebra.

This gives the first---``dynamical''---description of the Lie algebra
$\TGA(\CX,T)$.

\subsection{Lie algebras with root system $(\CX,T)$}
Definition of Lie algebra will be followed to Kac--Moody pattern but
with important changes.  First of all we define a \emph{local
algebra}.  Let $\phi\in C(\CX)$; we consider three types of uncountably
many generators: $X_{-1}(\phi)$, $X_0(\phi)$, $X_{+1}(\phi)$ where
$\phi$ runs over $C(\CX)$.  The list of relations is as follows:
\begin{equation}\label{eq3}
\begin{split}
[X_0(\phi),X_0(\psi)]&\mbox{}=0,\\
[X_0(\phi),X_{\pm1}(\psi)]&\mbox{}=\pm X_{\pm1}(K\phi\cdot\psi),\\
[X_{+1}(\phi),X_{-1}(\psi)]&{}=X_0(\phi\cdot\psi),
\end{split}
\end{equation}
where product ($\cdot$) is the product in associative algebra $C(\CX)$
and $K$ is a linear operator in $C(\CX)$ which is called Cartan
operator:
\begin{equation}\label{eq4}
(K\phi)(x)=2\phi(x)-\phi(Tx)-\phi(T^{-1}x).
\end{equation}
It is evident that Jacobi identity is true (if it makes sense) in the
local algebra
\[
\GA_{-1}\oplus\GA_0\oplus\GA_{+1},
\]
where $\GA_i=\Span\{X_i(\phi),\,\phi\in C(\CX)\}$, $i=-1,0,+1$.
($\GA_0\simeq C(\CX)$ is Cartan subalgebra.)  The further steps are the
same as in Kac--Moody theory \cite{K}.

We define free Lie algebra which is generated by the local algebra and
factorize it over the maximal ideal which has zero intersection with
$\GA_0$.  The resulting Lie algebra is denoted by $\GA(\CX,T)$.  The
fact is that this algebra is the same as $\TGA(\CX,T)$ of subsection 2.

We omit the verification that $\GA(\CX,T)$ is the graded Lie algebra
with the graded structure (as a linear space) as follow
\[
\bigoplus_{n\in\BZ}\GA_n
\]
and each $\GA_n\simeq C(\CX)$ (see \cite{V}).
So we can denote the
elements of $\GA_n$ as $X_n(\phi)$, $\phi\in C(\CX)$.

\begin{theorem}
The following formulas give the canonical isomorphism $\tau$ between
$\TGA(\CX,T)$ and dense part of $\GA(\CX,T)$:

\[
\begin{split}
\tau(\phi\ot U^n)&{}=\begin{cases}
X_{-n}(U^{-n}\phi), & (n>0) \\
X_0(\phi-U^{-1}\phi), & n=0,\quad \int_\CX\phi\,d\mu=0\\
X_n(\phi) & (n>0) \end{cases}\\
\tau(\mathbf{1}\ot U^0)&{}=X_0(\mathbf{1}).
\end{split}
\]
\end{theorem}

\begin{proof} The kernel of $\tau$ is $\mathbf{0}$.
Let us check that $\tau([a,b])=[\tau a,\tau b]$.  It is enough to test
monomials only.
\begin{multline*}
\bigl[\tau(\phi\ot U^0),\tau(\psi\ot U)\bigr]
=\bigl[\phi-U^{-1}\phi)\ot U^0,\psi\ot U\bigr]\\
=\psi\bigl(\phi-U^{-1}\phi-U(\phi-U^{-1}\phi)\bigr)\ot U
=(\psi\cdot K\phi)\ot U\\
=X_{+1}(K\phi\cdot\psi)
=\tau\bigl([\phi\ot U^0,\psi\ot U]\bigr);
\end{multline*}
\begin{multline*}
\bigl[X_0(\phi),X_{-1}(\psi)\bigr]
=\bigl[(\phi-U^{-1}\phi)\ot U^0, U^{-1}\psi\ot U^{-1}\bigr]\\
=\bigl(U^{-1}\psi(\phi-U^{-1}\phi-U^{-1}\phi+U^{-2}\phi)\bigr)\ot U^{-1}
=\bigl(U^{-1}(\psi(U\phi-2\phi+U^{-1}\phi))\bigr)\ot U^{-1}\\
=-U^{-1}(\psi\cdot K\phi)\ot U^{-1}
=-X_{-1}(K\phi\cdot\psi);
\end{multline*}
\begin{multline*}
\bigl[X_{+1}(\psi),X_{-1}(\gamma)\bigr]
=\bigl[\psi\ot U,U^{-1}\gamma\ot U^{-1}\bigr]\\
=(\psi\gamma-U^{-1}\gamma\cdot U^{-1}\psi)\ot U^0
=\bigl(\psi\gamma-U^{-1}(\psi\gamma)\bigr)\ot U^0
=X_0(\psi\gamma);
\end{multline*}
\[
\bigl[X_{+1}(\mathbf{1}),X_{-1}(\mathbf{1})\bigr]
=c X_0(\mathbf{1})=c\mathbf{1}.
\]
\end{proof}

\begin{remark} 1. We calculated the bracket only for elements of
local algebra, but this is enough because it generates all algebra.

\noindent 2. The $\tau$-image of $\TGA(\CX,T)$ is not all $\GA(\CX,T)$
but dense part of $\GA(\CX,T)$, for example, the set of functions
$\phi-U\phi$ is dense in $C_0(\CX)$ only, but $\GA(\CX,T)$ is the
extension of the image of $\TGA(\CX,T)$ and we can consider $\GA(\CX,T)$
as some kind of completion of $\TGA(\CX,T)$.
\noindent 3. Using isomorphisms $\tau$ we can rewrite
the definition of cocycle
of $\GA(\CX,T)$ ($\tau(\phi\ot U^{-n})=X_{-n}(U^{-n}\phi)$, $n>0$) so
\[
\alpha\bigl(X_n(\phi),X_m(\psi)\bigr)=
\begin{cases} 0, & \text{if $n+m\ne0$,}\\
n\int_{\CX}\phi\psi\,d\mu, & n+m=0. \end{cases}
\]
\end{remark}

Lie algebra $\GA(\CX,T)$ defined above is our main object.

\begin{proof}
If $n>0$ then $\alpha (X_n(\phi),X_{-n}(\psi))
=n\int_\CX\phi\cdot U^n(U^{-n}\psi)d\mu=n\int_\CX\phi\psi\,d\mu$.
If $n<0$ the $\alpha (X_{-n}(\phi),X_{n}(\psi))
=-n\int_\CX U^{-n}\phi\cdot U^{-n}\psi\,d\mu=-n\int_\CX\phi\psi\,d\mu$
($U$ is an unitary  operator).
\end{proof}

This consists with the initial formula
\[
\bigl[X_{+1}\bigl(\phi\bigr),
X_{-1}\bigl(({\phi})^{-1}\bigr)\bigr]=X_0(\mathbf{1})=\mathbf{1}\cdot c
\]
for $n=\pm1$.

Now we can rewrite the brackets for all monomials (not only for local part).
Assume $n,m>0$.
{\renewcommand{\theequation}{+,+}
\begin{equation}
\bigl[X_n(\phi),X_m(\psi)\bigr]
=\bigl[\phi\ot U^n, \psi\ot U^m\bigr]
=X_{n+m}(\phi\cdot U^n\psi-\psi\cdot U^m\phi);
\end{equation}
}
{\renewcommand{\theequation}{+,$-$}
\begin{equation}
\begin{split}
\bigl[X_n(\phi)&,X_{-m}\psi\bigr]
=\bigl[\phi\ot U^n, U^{-m}\psi\ot U^{-m}\bigr]\\
&{}=\bigl(\phi\cdot U^{n-m}\psi-U^{-m}(\phi\psi)\bigr)\ot U^{n-m}\\
&{}=\begin{cases}
X_{n-m}\bigl(\phi\cdot U^{n-m}\psi-U^{-m}(\phi\psi)\bigr), & n>m>0,\\
X_{0}\bigl((1-U^{-m})(1-U^{-1})^{-1}\phi\psi\bigr), & n=m,\\
X_{n-m}\bigl(U^{-n}\phi(\psi\cdot U^m\phi-U^{-n}\psi)\bigr), & 0<n<m;
\end{cases}
\end{split}
\end{equation}
}
{\renewcommand{\theequation}{$-$,$-$}
\begin{equation}
\begin{split}
\bigl[X_{-n}(\phi)&,X_{-m}\psi\bigr]
=\bigl[U^{-n}\phi\ot U^{-n},U^{-m}\psi\ot U^{-m}\bigr]\\
&=(U^{-n}\phi\cdot U^{-n-m}\psi-U^{-m}\psi\cdot U^{-n-m}\phi)\ot U^{-n-m}\\
&=X_{-n-m}(\psi\cdot U^m\phi-\phi\cdot U^n\psi)\\
&=-X_{-n-m}(\phi\cdot U^n\psi-\psi\cdot U^m\phi)
\end{split}
\end{equation}
}
{\renewcommand{\theequation}{0,+}
\begin{equation}
\begin{split}
\bigl[X_0(\phi),&X_n(\psi)\bigr]
=\bigl[(\phi-U^{-1}\phi)\ot U^0, \psi\ot U^n\bigr]\\
&=\bigl(\psi\cdot(\phi-U^{-1}\phi+U^{n-1}\phi-U^n\phi)\bigr)\ot U^n
=X_n(K_n\phi\cdot\psi),
\end{split}
\end{equation}
}
where $K_n=I-U^{-1}+U^{n-1}-U^n=(I-U^{-1})(I-U^n)$;
{\renewcommand{\theequation}{0,$-$}
\begin{equation}
\begin{split}
\bigl[X_0(\phi),&X_{-n}\psi\bigr]
=\bigl[(\phi-U^{-1}\phi)\ot U^0, U^{-n}\psi\ot U^{-n}\\
&=\bigl(U^{-n}\psi\cdot(\phi-U^{-1}\phi-U^{-n}\phi+U^{-n-1}\phi)\bigr)
\ot U^{-n}\\
&=X_{-n}((U^{-n}\psi)\cdot (U^n\phi-U^{n-1}\phi-\phi+U^{-1}\phi)
=-X_n(K_n\phi\cdot\psi)
\end{split}
\end{equation}
}
We can now observe that the formulas (+,+) and ($-$,$-$) are the same,
as well as (0,+) and (0,$-$).

\begin{theorem}
The formulas for the brackets of monomials in the subalgebra $\GA(\CX,T)$
are the following:

\noindent
\emph{1)} $\bigl[X_n(\phi),X_m(\psi)\bigr]
=\pm X_{n+m}(\phi\cdot U^{\pm n}\psi-\psi U^{\pm m}\phi)$,\\
where the sign is ``$+$'' if $n,m>0$ and ``$-$'' if $n,m<0$.

\noindent
\emph{2)} $\bigl[X_0(\phi),X_{\pm n}(\psi)\bigr]=\pm X_n(K_n\phi\cdot\psi)$.

\noindent
\emph{3)} $\displaystyle\bigl[X_n(\phi), X_m(\psi)\bigr]
=\begin{cases}
X_{n+m}\bigl(U^m\psi(\phi\cdot U^n\psi-U^m\phi)\bigr), &
{\scriptstyle  \begin{cases}\scriptstyle m<0,\\
\scriptstyle n+m>0\end{cases}}\\
X_0\bigl(\frac{1-U^{-m}}{1-U^{-1}}
(\phi\psi)\bigr)+n\int_\CX\phi\psi\,d\mu &
{\scriptstyle \begin{cases}\scriptstyle n+m=0,\\
 \scriptstyle n\ne0\end{cases}}\\
X_{n+m}\bigl(U^{-n}\phi(\psi\cdot U^{-m}\phi-U^{-n}\psi)\bigr), &
{\scriptstyle \begin{cases}\scriptstyle m<0,\\ \scriptstyle n+m<0\end{cases}}
\end{cases}
$

\noindent
\emph{4)} $\bigl[X_0(\phi),X_0(\psi)\bigr]=0$.\qed
\end{theorem}

Lie algebra $\GA(\CX,T)$ does not associate with associative
algebra; cocycle $\alpha$ has nothing to do with associated
crossproduct of subsection 1.1.  The role of central extension
is very important.

We defined Lie algebra $\GA(\CX,T)$ ($\simeq\TGA(\CX,T)$) in a new
terms, compare with (3).  This manner gives us the formulas for
local part (1--2)
which are similar to classical ones (Cartan simple algebras and
Kac--Moody algebras).  But the formulas for general monomials
are more complicated than dynamical (see subsection 1.1) description.

\section{General Lie algebras with continuous root systems
and new examples of $\GA(\CX,T)$}

\subsection{General definition}
We recall \cite{VS1,VS2,V} the  definition of graded Lie algebras
with continuous root system.

Suppose $\CH$ is a commutative associative Lie $\BC$-algebra with
unity (Cartan subalgebra) and $K:\CH\hookleftarrow$ is a linear operator
(Cartan operator).  The local algebra $[K]$ is, as a linear space, a
 direct sum
\[
\CH_{-1}\oplus\CH_0\oplus\CH_{+1},\qquad \CH_{i}\simeq\CH,\quad i=0,\pm1
\]
with brackets:
\[
\begin{split}
X_i(\phi)\in\CH_i,\quad i=0,&\pm1,\quad\phi,\psi\in\CH\\
\bigl[X_0(\phi),X_0(\psi)\bigr]&{}=0,\\
\bigl[X_0(\phi),X_{\pm1}(\psi)\bigr]&{}=\pm X_{\pm1}(K\phi\cdot\psi)\\
\bigl[X_{+1}(\phi),X_{-1}(\psi)\bigr]&{}=X_0(\phi\cdot\psi)
\end{split}
\]
The local algebra $\CH_{-1}\oplus\CH_0\oplus\CH_{+1}$ generates
graded Lie algebra $\GA(\CH,K)$.
in the same spirit as in Subsection~1.2 (and as in the theory of
LKM-algebras).
Then we obtain $\GA(\CX,T)$ from Section~1.

The spectrum of commutative algebra $\CH$ (if it exists) is root system
of $\GA(\CH,K)$
by definition (see \cite{V}), more exactly the set of
simple positive roots. But it could be no spectra (say, $\CH$ is the
algebra of rational functions) so we have Lie algebras without simple
roots but with Cartan operator.

The condition of constant or polynomial growth of the dimension (in
an appropriate sense) puts essential restriction on the operator $K$.

\begin{remark}
Let $E\subset\CH$ is an invariant under $K$ subalgebra of $\CH$.
Then $\GA(E,K)$ is Lie subalgebra of $\GA(\CH,K)$.  In particular, if
$E_1\subset E_2\subset\dots$, $\cup_i E_i=\CH$, is a sequence of
$K$-invariant subalgebras of $\CH$ then
$\GA(\CH,K)=\cup_{i=1}^\infty\GA(E_i,K)$.
\end{remark}

\subsection{New examples of algebras of type $\GA(\CX,T)$}
The first nontrivial example of algebras of type $\GA(\CX,T)$ was
so called sine-algebra.  We will not consider it because it was done
before from different point of view (see \cite{VS2,FFZ,GKL,V}).
It was defined independently in \cite{VS2} and \cite{FFZ}. We give
now general example of similar type.

Let $(\CX,T,\mu)$ be an ergodic system with discrete spectrum. It
means that operator $U=U_T$ has spectral decomposition
\[
Uf=\sum_{\lambda}\lambda f_\lambda\chi_\lambda,\text{ where }
f=\sum_\lambda f_\lambda\chi_\lambda,
\]
sum is over eigenvalues of $U$ ($\lambda\in\BT^1$), and $\chi_\lambda$
is the eigenfunction corresponding to $\lambda$.  It is well-known
(von Neumann theorem) that such system can be realized on the
compact abelian group $G=\CX$ with Haar measure $m=\mu$ and $T$ is
translation on some element $g_0\in G$, then $\chi_\lambda$ is a
character of $G$ ($\chi_\lambda\in G^\wedge$) and
\[
\lambda=\chi_\lambda(g_0)\in\BT^1.
\]
Sine-algebra corresponds to the case $\CX=\BT^1$ and $T$ is translation on
irrational number $\theta\in S^1=\BR/BZ$.  In \cite{GKL} was considered
also the case $\CX=\BT^d\ni\theta$.

The case
\[
\CX=Z_p,\qquad Tx=x+1,
\]
where $Z_p$ is additive group of $p$-adic integers, $p$ is a prime, and
$T$ is adding of unity, is more interesting from our point of view.
The measure $\mu=m$ is Haar (additive) measure on $Z_p$.

The group of characters $Z_p^\wedge=Q_{p\infty}$ is a group of all roots
of unity of the degree $p^n$, $n\in\BN$, and the characters
$\mu\in Q_{p\infty}$ (as function on $Z_p$ with values in $\BT^1$)
are eigenfunctions of operator $U=U_T$.

We give the description of the general case when operator $U=U_T$
has discrete spectrum. The specific property of algebra $\GA(\CX,T)$ in
this case is existence of the natural \emph{linear basis} in $\GA(\CX,T)$.

Suppose $G$ is abelian (additive) compact group and $G^\wedge$ is a
countable group of the (multiplicative) characters on $G$.   We fix
the element $\lambda\in G$ with dense set of powers: $\Cl\{\lambda^n,
n\in\BZ\}=G$.  Then $T_\lambda g=Tg=g+\lambda$, $U_T=U$,
$(Uf)(g)=f(g+\lambda)$, $f\in C(G)$.  Each character $\chi\in G^\wedge$
is an eigenfunction of $U$ with eigenvalue $\chi(\lambda)\in\BT^1$.

\begin{theorem}
Linear basis  in the Lie algebra $\GA(G,T)$ is the set
$\{Y_{\chi,n}:\chi\in G^\wedge, n\in\BZ\}$ with the following brackets
\setcounter{equation}{4}
\begin{equation}\label{eq5}
\bigl[Y_{\chi,n},Y_{\chi_1,n_1}\bigr]
=\bigl(\chi_1(\lambda)^n-\chi(\lambda)^{n_1}\bigr) Y_{\chi\chi_1,n+n_1}
  +\delta_{n+n_1}\tilde{\delta}_{\chi\chi_1}\cdot n\cdot c,
\end{equation}
where $\delta_n=\begin{cases} 1, & n=0\\ 0, &n\ne0\end{cases}$,
$\tilde{\delta}_\chi=\begin{cases} 1, & \chi=\mathbf{1}\\
   0, &\chi\ne\mathbf{1}\end{cases}$.

Algebra $\GA(G,T)$ is $\BZ\times G^\wedge$-graded algebra; the
subalgebra $\{c1:c\in\BC\}$
in Cartan subalgebra $\GA_0=C(G)$ is the
center of $\GA(G,T)$.
\end{theorem}

\begin{proof}
Assume $Y_{\chi,n}=\chi\ot U^n$ as an element of $\GA_n$, where $\chi$
is the character of $G$ as a function $G\to\BT^1$. It is easy to check
that the brackets (see formula \eqref{eq1} in Section~1) give us
formula \eqref{eq5}.  Note that the center is not direct summand, so
we have
\[
\bigl[Y_{\chi,n},Y_{\chi^{-1},-n}\bigr]=n\cdot c.\qed
\]
\renewcommand{\qed}{}
\end{proof}

This is the third description of our algebra $\GA(\CX,T)$ with linear basis;
this description is valid for discrete spectrum only.

The group $G$ is the set of simple roots for $\GA(G,T)$ and ``Dynkin''
diagram is the set of arrows $G\ni g\to g+\lambda\in G$.  In opposite
to Kac--Moody case our algebras $\GA(\CX,T)$ have no imaginary roots.

\subsection{The case of $p$-adic integers}
Return back to the case $G=Z_p$, $Tx=x+1$.  In this case $\GA(Z_p,T)$ is
$\BZ\times Q_{p\infty}$-graded algebra. Cartan subalgebra is space $C(Z_p)$.
Consider finite dimensional subspaces $L_n\subset C(Z_p)$ of functions
depending on the points of the quotient $Z_p\to\BZ/p^n$; it is clear
that subspace $L_n$ is $U_T$-invariant, so
$\CL_n=\bigoplus_{m\in\BZ}L_n\ot U^m$

is a subalgebra of $\GA$.

\begin{theorem}
The Lie algebra $\CL_n$ is canonically isomorphic to the algebra
$\CA^{(1)}_{p^n}$.  Consequently, $\GA(Z_p,T)$ is (completion) of the
inductive limit of Kac--Moody Lie algebras\footnote{More exactly,
inductive limit contains only linear combinations of monomials of type
$\phi\ot U^n$, where $\phi$ are cylindric functions; so it is enough to
extend the set of $\phi$ onto arbitrary continuous functions which
means to make a completion.} \[ \GA(Z_p,T)\supset \limind_n\CL_n.\qed
\]
\renewcommand{\qed}{}
\end{theorem}

\begin{remark} It is possible to define Weyl group $W$ for this algebra.
Group $W$ contains the group of permutations of the coordinates in $Z_p$.

It is very instructive to study the link between theory of Kac--Moody
affine algebras and our theory looking on this example.
\end{remark}

\end{document}